\documentclass[12pt]{amsart}

\usepackage[textheight=23.5cm, left=1.2cm, right=1.2cm]{geometry}

\usepackage[english]{babel}
\usepackage[T1]{fontenc}
\usepackage[utf8]{inputenc}

\theoremstyle{plain}
    \newtheorem{theorem}{Theorem}[section]
    \newtheorem{lemma}[theorem]{Lemma}

\theoremstyle{definition}

\usepackage[
    draft = false,
    unicode = true,
    colorlinks = true,
    allcolors = blue,
    hyperfootnotes = true,
    citecolor = red
]{hyperref}
\usepackage{amsfonts,amsmath,amssymb,amscd,amsthm,url}
\usepackage[inline]{enumitem}
\usepackage{graphicx,epsf,afterpage}
\usepackage{soul}
\usepackage{multicol}
\usepackage{array}
\usepackage{braket}
\usepackage{epigraph}
\usepackage{rotating}

\usepackage[
backend=biber,
style=numeric-comp,
sorting=none,
giveninits=true
]{biblatex}
\usepackage{csquotes}

\usepackage{xcolor}

\usepackage{ulem}

\usepackage{tikz}
\usetikzlibrary{positioning, calc, intersections, through, arrows, matrix, chains, math}
\usetikzlibrary{shadows}
\usetikzlibrary{graphs}
\usetikzlibrary{graphs.standard}
\usetikzlibrary{patterns}
\usetikzlibrary{decorations.pathreplacing,calligraphy,backgrounds}

\newcommand\norm[1]{\ensuremath{\left\lVert#1\right\rVert}}

\renewcommand{\Pr}{\mathrm{P}}

\DeclareMathOperator{\Expect}{\mathbb{E}}

\newcommand{\Lcal}{\mathcal{L}}
\newcommand{\Hcal}{\mathcal{H}}
\newcommand{\Tcal}{\mathcal{T}}

\newcommand{\Xcal}{\mathcal{X}}

\newcommand{\R}{\ensuremath{\mathbb{R}}}

\newcommand{\N}{\ensuremath{\mathbb{N}}}

\renewcommand{\geq}{\geqslant}

\renewcommand{\leq}{\leqslant}

\newcounter{mcnt}

\newcounter{wordcnt}

\begin{document}

\title[The SLLN for random semigroups with unbounded generators on uniformly smooth Banach spaces]
{The Strong Law of Large Numbers \\ for random semigroups with unbounded generators \\ on uniformly smooth Banach spaces}

\author{S.\,V. Dzhenzher and V.\,Zh. Sakbaev}

\begin{abstract}
     We consider random linear unbounded operators on a Banach space $\mathcal{X}$.
     For example, such random operators may be random quantum channels.
     The Law of Large Numbers is known when $\mathcal{X}$ is a Hilbert space, in the form of the usual Law of Large Numbers for random operators, and in some other particular cases.
     Instead of the sum of i.i.d. variables, there may be considered the composition of random semigroups $e^{A_i t/n}$.
     We obtain the Strong Law of Large Numbers in Strong Operator Topology for random semigroups of unbounded linear operators on a uniformly smooth Banach space.
\end{abstract}

\thanks{\hspace{-4mm}S.\,V. Dzhenzher: sdjenjer@yandex.ru. orcid: 0009-0008-3513-4312
\\
V.\,Zh. Sakbaev: fumi2003@mail.ru. orcid: 0000-0001-8349-1738
\\
V.\,Zh. Sakbaev: Keldysh Institute of Applied Mathematics of Russian Academy of Sciences 125047, Moscow, Russia
\\
Both authors: Moscow Institute of Physics and Technology 141701, Dolgoprudny, Russia}

\newcommand{\LO}{\mathcal{LO}}

\maketitle
\thispagestyle{empty}


\vspace{5mm}


\section{Introduction}

The theory of compositions of random linear operators is constructed and developed in lots of publications \cite{Mehta, Furstenberg, Tutubalin, Skorokhod}. 
A random quantum dynamics can be considered in two ways: either as a composition of random linear operators in a Hilbert space \cite{Kempe, OSS19}, or as a composition of random quantum channels \cite{Aharonov-D-Z, Holevo-qprobstat, GOSS-2022, Pathirana-Schenker, Dhamapurkar-Dahlsten}.
See some historical exposition in \cite{DzhenzherSakbaev25-SLLN, Dzhenzher25}.

The Law of Large Numbers (LLN) and the Central Limit Theorem (CLT) for random variables acting in a Banach space are well studied \cite{Ledoux-Talagrand-1991-probability}.
Distributions of products of random linear operators have been studied in \cite{Oseledets, Tutubalin, Skorokhod, Pathirana-Schenker}.
Different versions of limit theorems of LLN and CLT type for compositions of i.i.d. random operators were obtained in \cite{Berger, Watkins-1984, Watkins-1985, DzhenzherSakbaev25-SLLN, Dzhenzher25, SSSh23}.

However, a consideration of quantum dynamics of an open quantum system needs to study the random operators in the Banach space of quantum states \cite{Holevo-qprobstat, GOSS-2022, Kalmetev-Orlov-Sakbaev-2023, Pechen, Sahu}. Nevertheless, the case of random operators acting in a Banach space is studied significantly less. This paper aims to fill this gap for the Strong Law of Large Numbers (SLLN) for random semigroups of unbounded linear operators in a uniformly smooth Banach space.
In particular, in \cite{DzhenzherSakbaev25-SLLN}, the SLLN was obtained for them in the spirit of Chernoff's Theorem \cite{Chernoff-1968}; it was obtained in the Weak Operator Topology for all Banach spaces, and in the Strong Operator Topology (SOT) for uniformly smooth Banach spaces.
It is only natural to consider random semigroups with unbounded generators, following the classics \cite{Kato, Engel-Nagel-sec2, Engel-Nagel-sec3}.

The structure of the article is as follows.
In~\S\ref{s:basic}, basic notions and the main result (Theorem~\ref{t:slln-unbound}) are given.
In~\S\ref{s:proof}, the proof of Theorem~\ref{t:slln-unbound} is given.
In~\S\ref{s:concl}, the discussion of the obtained results is given.

\section{Basic notions and results}\label{s:basic}

Recall that a Banach space \(\mathcal{X}\) is said to be \textbf{\emph{uniformly smooth}} if the limit
\[
    \lim_{t\to 0} \frac{\norm{x + ty} - \norm{x}}{t}
\]
exists for all \(x,y\) from the unit sphere in \(\mathcal{X}\).
For example, for \(1 < p < \infty\) the spaces \(L_p\) of functions whose $p$-th power is integrable, the Sobolev spaces \(W^{1,p}(0,1)\), and the Schatten classes \(\mathcal{T}_p\) are uniformly smooth; and the spaces \(L_1\), \(C[\,0,1\,]\), and \(c_0 \subset \ell_\infty\) are not uniformly smooth.
By \cite[Theorem~4.24]{Pisier-2016-martingales}, all uniformly smooth Banach spaces are exactly \textbf{\emph{$p$-smooth}} Banach spaces for \(p \in (1,2\,]\), which means that there is $C >0$ such that for any \(x,y\in \mathcal{X}\) we have
\[
    \norm{x+y}^p + \norm{x-y}^p \leq 2\norm{x}^p + C\norm{y}^p.
\]
It is well known that all Hilbert spaces are $2$-smooth.
Also, for \(1 < p < 2\) the spaces \(L_p\) are $p$-smooth, and for \(2 \leq p < \infty\) the spaces \(L_p\) are $2$-smooth.
See details in \cite{Ledoux-Talagrand-1991-probability, DzhenzherSakbaev25-SLLN}.

Let \((\Omega,\sigma_\Omega,\Pr)\) be a probability space; that means, $\Omega$ is a set of sample points, $\sigma_\Omega$ is a $\sigma$-algebra over $\Omega$, and \(\Pr\colon\sigma_\Omega\to [\,0,1\,]\) is a $\sigma$-additive measure.

For a Banach space $\Xcal$ denote by \(\Lcal(\Xcal)\) the space of linear bounded operators on \(\Xcal\).
Now, denote by \(\LO(\Xcal) \supset \Lcal(\Xcal)\) the space of linear densely defined but not necessary bounded operators on \(\Xcal\).
Stricly speaking, \(A\in \LO(\Xcal)\) is an operator \(A\colon D(A) \to \Xcal\) defined on the dense domain \(D(A) \subset \Xcal\).

In this section we consider linear but not necessarily bounded random operators \(\Omega\to \LO(\Xcal)\) on a uniformly smooth Banach space $\Xcal$.
We mean that the domain \(D(A)\subset\Xcal\) of a random operator $A$ is independent of $\omega\in\Omega$.
By measurability of a random operator $(A,D(A))$, we mean \emph{strong measurability}, which means that the random element \(Ax\) is Bochner measurable for any \(x\in D(A)\).

The expected value \(\Expect\xi\) of a random element \(\xi\) is defined as its Bochner integral; in particular, this works for \(\xi\) being a linear bounded operator on a Banach space.
The expected value \(\Expect L\) of a random operator $(L,D(L))$ is defined pointwise as the Bochner integral: that is, for any \(x\in D(L)\)
\[
    (\Expect L)x := \Expect (Lx).
\]

Let $(L_0, D(L_0))\in\LO(\Xcal)$ be a non-random operator such that $-L_0$ is closable and dissipative; the latter means that
\[
    \norm{(\lambda + L_0)x} \geq \lambda\norm{x}
\]
for all $\lambda>0$ and \(x\in D(L_0)\).
Suppose that the closure of $-L_0$ generates a strongly continuous semigroup $e^{-L_0t}$ satisfying
\[
    \norm{e^{-L_0t}} \leq Me^{\beta t}
\]
for some $M\geq 1$ and \(\beta\in\R\).

For $C>0$ denote by \(B_C(L_0) \subset \LO(\Xcal)\) the class of linear operators $L$ with the domain \(D(L) = D(L_0)\) and such that
\[
    \norm{L-L_0} \leq C;
\]
in particular, this means that $L-L_0$ is defined on the whole $\Xcal$, and hence lies in \(\Lcal(\Xcal)\).

We now consider random operators \(L\colon\Omega\to B_C(L_0)\).
By \cite[Theorem~2.1]{Kato} (see also \cite[Bounded Perturbation Theorem~1.3]{Engel-Nagel-sec3}), they generate random strongly continuous semigroups $e^{-Lt}$ satisfying
\[
    \norm{e^{-Lt}} \leq Me^{(\beta + MC)t}.
\]
We denote \(\gamma:=\max(\beta+MC,1)\) and denote by \(\Lcal_{M,\gamma} \subset B_C(L_0)^\Omega\) the set of random operators $L$ such that semigroups $e^{-Lt}$ are strongly measurable and such that \(\Expect L = L_0\).

\begin{theorem}[Chernoff-type SLLN for i.i.d. random semigroups; proved in \S\ref{s:proof} below]\label{t:slln-unbound}
    Let $\Xcal$ be a uniformly smooth Banach space.
    Let \(L_1,L_2,\ldots \subset \Lcal_{M,\gamma}\) be i.i.d. random operators; this means that there exists \(\Expect L_i = L_0 \in \LO(\Xcal)\) such that \(L_i-L_0\colon\Omega\to\Lcal(\Xcal)\), and that semigroups \(e^{-L_it}\) are strongly measurable and satisfy
    \[
        \norm{e^{-L_it}} \leq Me^{\gamma t}.
    \]
    Then \(e^{-L_1t/n}\ldots e^{-L_nt/n}\) converges a.s. in SOT to \(e^{-L_0 t }\) uniformly for $t$ in any segment; that is, for any \(x\in\Xcal\) and $T>0$ almost surely
    \[
        \lim_{n\to\infty}\sup_{t\in[\,0,T\,]} \norm{(e^{-L_1t/n}\ldots e^{-L_nt/n} - e^{-L_0 t})x} = 0.
    \]
\end{theorem}

Recall that a \emph{quantum state} on a separable Hilbert space \(\Hcal\) is a trace-class operator from the $1$-st Schatten class \(\Tcal_1(\Hcal)\); sometimes it is additionally required for a quantum state to be positive and with unit trace.
Recall that a \emph{quantum channel} on a separable Hilbert space \(\Hcal\) is a completely positive trace-preserving linear map \(\Tcal_1(\Hcal)\to\Tcal_1(\Hcal)\).
If \(\Hcal\) is $d$-dimensional, then \(\Tcal_1(\Hcal)\) is just the space of \emph{density matrices}.
For \(0 < \lambda< 1 + \frac{1}{d^2-1}\), denote by \(D_\lambda\colon\Tcal_1(\Hcal)\to\Tcal_1(\Hcal)\) the \emph{$d$-dimensional quantum depolarizing channel} given by
\[
    D_\lambda \rho := (1-\lambda)\rho + \frac{\lambda}{d}I,
\]
where \(I/d\) is the \emph{maximally mixed state}.
As we could see, with \(\lambda\to +0\) the channel \(D_\lambda\) could be approximated as \(e^{\lambda A/d}\), where $A$ is given by \(A\rho = I/d\) for all \(\rho\).
That is why the following is a version of the previous theorem in relation to quantum channels.

\begin{theorem}[Chernoff-type SLLN for quantum channels]\label{t:slln-quch}
    Let \(\{\xi_n\}\) be i.i.d. random variables on the interval \((0, 1 + \frac{1}{d^2-1})\).
    Then the random channels
    \[
        D_{\xi_1/n}\ldots D_{\xi_n/n}
    \]
    converge almost surely in SOT to \(D_{1-e^{-\Expect\xi_1}}\).
\end{theorem}

\begin{proof}
    By definition of a quantum depolarizing channel, for any quantum state \(\rho\) we have
    \[
        D_{\xi_1/n}\ldots D_{\xi_n/n} \rho =
        \rho\prod_{k=1}^n \left(1-\frac{\xi_k}{n}\right) + \frac{I}{d}\sum_{k=1}^n \frac{\xi_k}{n}\prod_{i=1}^{k-1} \left(1-\frac{\xi_i}{n}\right).
    \]
    Now
    \[
        \prod_{k=1}^n \left(1-\frac{\xi_k}{n}\right) \xrightarrow[n\to\infty]{a.s.} e^{-\Expect\xi_1}
    \]
    and
    \[
        \sum_{k=1}^n \frac{\xi_k}{n}\prod_{i=1}^{k-1} \left(1-\frac{\xi_i}{n}\right) \xrightarrow[n\to\infty]{a.s.} 1-e^{-\Expect\xi_1}
    \]
    imply the result.
\end{proof}

\section{Proof of Theorem~\ref{t:slln-unbound}}\label{s:proof}

Denote by \(I\in\Lcal(\Xcal)\) the identity operator.

The following theorem is the analogue of the well-known Chernoff Product Formula \cite{Chernoff-1968}.

\begin{theorem}[Chernoff Product Formula for i.i.d. random operators]\label{t:chernoff}
    Let $\Xcal$ be a Banach space.
    Let \(L_1,L_2,\ldots \subset \Lcal_{M,\gamma}\) be i.i.d. random operators; this means that there exists \(\Expect L_i = L_0 \in \LO(\Xcal)\) such that \(L_i-L_0\colon\Omega\to\Lcal(\Xcal)\), and that semigroups \(e^{-L_it}\) are strongly measurable and satisfy
    \[
        \norm{e^{-L_it}} \leq Me^{\gamma t}.
    \]
    Then \(\Expect e^{-L_1t/n}\ldots e^{-L_nt/n}\) converges in SOT to \(e^{-L_0 t }\) uniformly for $t$ in any segment; that is, for any \(x\in\Xcal\) and $T>0$
    \[
        \lim_{n\to\infty}\sup_{t\in[\,0,T\,]} \norm{(\Expect e^{-L_1t/n}\ldots e^{-L_nt/n} - e^{-L_0 t})x} = 0.
    \]
\end{theorem}

\begin{proof}
    Denote the~mapping $F\colon [\,0, +\infty)\to\Lcal(\Xcal)$ by the rule $F(t) := \Expect e^{-L_1t}$.
    By the Chernoff Product Formula \cite[Corollary~5.3]{Engel-Nagel-sec3}, we have
    \[
        \lim\limits_{n\to\infty} \sup\limits_{t \in [\,0,T\,]} \norm{\left(e^{F'(0)t} - F\left(t/n\right)^n\right)x} = 0
    \]
    if
    (a) \(F(0) = I\),
    (b) \(\norm{F(t)^n} \leqslant Ke^{ant}\) for some $a \in\R$ and $K>1$, and
    (c) for any $y\in D\subset\Xcal$ there is \(F'(0)y := \lim\limits_{t \to +0} \dfrac{F(t)y - y}{t}\), where \(D\) and \((\lambda-F'(0))D\) are dense subspaces in \(\Xcal\) for some \(\lambda>a\).
    Here~(a) is obvious,
    (b) holds for $K=M$ and $a=\gamma$,
    and~(c) holds for $F'(0) = -L_0$ and $D=D(L_0)$ by \cite[Lumer--Phillips Theorem~3.15]{Engel-Nagel-sec2}.
    It remains to notice that the Chernoff Product Formula gives exactly the required equality since
    \(\Expect e^{-L_1t/n} \ldots e^{-L_nt/n} = F(t/n)^n\).
\end{proof}

Let \(L,L_1,\ldots,L_n \subset \Lcal_{M,\gamma}\) be i.i.d. operators.
Following \cite{DzhenzherSakbaev25-SLLN}, for $i \in \N$ and $s \geq 0$ denote
\[
    \Delta_i(s) := e^{-L_is} - \Expect e^{-Ls}.
\]
For integer $0 \leq k \leq n$ denote \([n]:= \{1,\ldots,n\}\), and denote by \(\binom{[n]}{k}\) the family of $k$-element subsets of \([n]\). 
For \(\{i_1 < \ldots < i_k\} \in \binom{[n]}{k}\) and \(s>0\) denote
\[
    F_{n,\{i_1, \ldots, i_k\}}(s) := (\Expect e^{-Ls})^{i_1-1} \Delta_{i_1}(s) (\Expect e^{-Ls})^{i_2-i_1-1} \ldots \Delta_{i_k}(s) (\Expect e^{-Ls})^{n-i_k},
\]
meaning
\[
    F_{n,\varnothing}(s) = (\Expect e^{-Ls})^n.
\]

\begin{lemma}\label{l:est-norm}
    Let $\Xcal$ be a Banach space.
    Let $L_1, \ldots, L_n \subset\Lcal_{M,\gamma}$ be i.i.d. random operators; this means that there exists \(\Expect L_i = L_0 \in \LO(\Xcal)\) such that \(L_i-L_0\colon\Omega\to\Lcal(\Xcal)\), and that semigroups \(e^{-L_it}\) are strongly measurable and satisfy
    \[
        \norm{e^{-L_it}} \leq Me^{\gamma t}.
    \]
    Then for each $n\in \N$ and \(P \in \binom{[n]}{k}\)
    \[
        \norm{F_{n,P}(s)} \leq (2\gamma s)^k e^{n\gamma s}.
    \]
\end{lemma}

\begin{proof}
    We have \(\Expect\norm{e^{-Ls}} \leq Me^{\gamma s}\).
    Hence by \cite[Theorem~II.4.ii]{DiestelUhl}
    \begin{equation}\label{eq:e-norm}
        \norm{\Expect e^{-Ls}} \leq Me^{\gamma s}.
    \end{equation}
    By \cite[Corollary~4.6.4 (the mean value theorem)]{BogachevSmolyanov}
    \[
        \norm{e^{-Ls} - I} \leq \gamma se^{\gamma s}.
    \]
    By~\eqref{eq:e-norm} and again by \cite[Corollary~4.6.4 (the mean value theorem)]{BogachevSmolyanov}
    \[
        \norm{\Expect e^{-Ls} - I} \leq \gamma se^{\gamma s}.
    \]
    Hence by the triangle inequality
    \begin{equation}\label{eq:delta-norm}
        \norm{\Delta_i(s)} \leq 2\gamma s e^{\gamma s}.
    \end{equation}
    Now the lemma follows as the binomial combinations of inequalities~(\ref{eq:e-norm}, \ref{eq:delta-norm}).
\end{proof}

Now we are ready for the proof of Theorem~\ref{t:slln-unbound}, which is analogous to the proof of \cite[Theorem~2.4]{DzhenzherSakbaev25-SLLN}.
In particular, it uses the following well-known result \cite[Theorem~4.52]{Pisier-2016-martingales}.

\begin{theorem}[Burkholder-type]\label{t:burkholder}
    Let $\mathcal{X}$ be a $p$-smooth Banach space.
    Let \(\mu_n = \sum_{k=1}^n d_k\) be a martingale.
    Then for any \(r \geq 1\) there is some \(C = C(\mathcal{X},r,p)>0\) such that
    \[
        \Expect\norm{\mu_n}^r \leq C\Expect\left(\sum_{k=1}^n \norm{d_k}^p\right)^{\frac{r}{p}}.
    \]
\end{theorem}

\begin{proof}[Proof of Theorem~\ref{t:slln-unbound}]
    For a fixed \(x\in\Xcal\) denote
    \[
        \mu_n(t) := \left(e^{-L_1t/n} \ldots e^{-L_nt/n} - \Expect e^{-L_1t/n} \ldots e^{-L_nt/n}\right)x =
        \sum_{k=1}^n\sum_{P\in\binom{[n]}{k}} F_{n,P}(t/n)x.
    \]
    By Chernoff Product Formula~\ref{t:chernoff} it is sufficient to prove that for any \(T>0\) a.s.
    \[
        \lim\limits_{n\to\infty} \sup\limits_{t \in [\,0,T\,]} \norm{\mu_n(t)} = 0.
    \]
    For \(k \leq n\) denote
    \[
        d_{n,k}(t) := \sum_{\substack{P \subset [n] \\ \max P = k}} F_{n,P}(t/n)x.
    \]
    Up to \(\norm{x}\) by Lemma~\ref{l:est-norm} we have
    \begin{equation}\label{eq:mart-dif}
        \norm{d_{n,k}(t)} \leq \sum_{j=1}^k \binom{k-1}{j-1} \left(\frac{2\gamma t}{n}\right)^j e^{\gamma t} =
        \frac{2\gamma t}{n} e^{\gamma t} \left(1+\frac{2\gamma t}{n}\right)^{k-1} \leq \frac{2\gamma t}{n}e^{3\gamma t}.
    \end{equation}
    Note that for fixed $t>0$ and \(n\in\N\) the sequence of
    \[
        \mu_m(t) = \sum_{k=1}^m d_{n,k}(t),
        \quad\text{where}\quad m\leq n,
    \]
    is a martingale w.r.t. \(\mathcal{F}_m = \sigma(L_1,\ldots,L_m)\), since
    \begin{multline*}
        \Expect^{\mathcal{F}_{m-1}} F_{n,\{i_1 < \ldots < i_j=m\}}(s) =
        (\Expect e^{-Ls})^{i_1-1} \Delta_{i_1}(s) (\Expect e^{-Ls})^{i_2-i_1-1} \ldots \\ \ldots \Delta_{i_{j-1}}(s) (\Expect e^{-Ls})^{i_j-i_{j-1}-1} \underbrace{\Expect \Delta_{i_j}(s)}_{=0} (\Expect e^{-Ls})^{n-i_j} = 0.
    \end{multline*}
    Since $\Xcal$ is uniformly smooth, it is $p$-smooth for some $p\in(1,2\,]$.
    Then by Burkholder-type Theorem~\ref{t:burkholder} for \(r:=\frac{2p}{p-1}\) and by inequality~\eqref{eq:mart-dif} up to the unimportant constant we have
    \[
        \Expect\norm{\mu_n(t)}^r \leq \left(n \left(\frac{2\gamma t}{n}\right)^pe^{3p\gamma t}\right)^{\frac{r}{p}} = \frac{1}{n^2} (2\gamma t)^r e^{3r\gamma t}.
    \]
    Thus by the Markov's inequality (up to the constant from Theorem~\ref{t:burkholder})
    \[
        \Pr\left\{ \norm{\mu_n(t)} > \varepsilon \right\} \leq \frac{\norm{x}^r}{\varepsilon^r} \frac{1}{n^2} (2\gamma t)^r e^{3r\gamma t}.
    \]
    Finally by the continuity from below of the probability we obtain
    \[
        \Pr\left\{\sup\limits_{t \in [\,0,T\,]} \norm{\mu_n(t)} > \varepsilon \right\} = O\left(n^{-2}\right).
    \]
    Now the desired result follows by the Borel--Cantelli lemma.
\end{proof}

\section{Discussion and conclusion}\label{s:concl}

As it was already said, the field of limit theorems for random operators is developed a lot.
There are several direction for the progress in this area.

For example, one could obtain the result, analogous to Theorem~\ref{t:slln-unbound}, but in the spirit of the Central Limit Theorem as in \cite{Dzhenzher25}.
This means that the whole difference \(e^{-L_1t/n}\ldots e^{-L_nt/n} - e^{-L_0 t}\) should be multiplied by \(\sqrt{n}\), and the convergence to a random Gaussian operator should be obtained.
If additionally the supremum on $t$ is written (as in Theorem~\ref{t:slln-unbound}), then possibly some technique using Brownian bridges is needed. This question is interesting even when the operators act on a Hilbert space.

Another possible way of development of these results is to obtain sharp restrictions on a Banach space.
It is well known that the geometry of a Banach space is tightly connected with the SLLN-type results in this space.
However, as it was shown in \cite{DzhenzherSakbaev25-SLLN}, for bounded random operators the uniform smoothness on a Banach space is not necessary for the SLLN-type results.
In particular, the SLLN-type result holds for bounded random operators on \(\ell_1\).
So, it would be intesting to weaken the assumption on uniform smoothness in Theorem~\ref{t:slln-unbound}, and to find counterexamples for other assumptions.

Finally, there is a more physical way to develop this field, in the spirit of Theorem~\ref{t:slln-quch}.
This theorem was easily proved without usage of any additional theory, since we considered only the particular channels on a finite-dimentional Hilbert space. It is interesting to obtain analogous results for other quantum channels, or under more general assumptions.

\printbibliography

\end{document}